\documentclass[twocolumn]{autart}    

\usepackage{graphicx}          
\usepackage{amssymb,amsmath}
\usepackage[norelsize,algo2e,linesnumbered,ruled,vlined]{algorithm2e}
\SetArgSty{textup}
\usepackage{subcaption}
\usepackage{url}

\begin{document}

\begin{frontmatter}

\title{Exploiting symmetries in active set enumeration for constrained linear-quadratic optimal control\thanksref{footnoteinfo}} 

\thanks[footnoteinfo]{Supports by the Alexander von Humboldt Foundation research group linkage cooperation program and the European Union’s Horizon Europe under grant no. 101079342 (Fostering Opportunities Towards Slovak Excellence in Advanced Control for Smart Industries) are gratefully acknowledged. M. Kvasnica gratefully acknowledges the contributions of the Slovak Research and Development Agency under the project APVV-20-0261 and the Scientific Grant Agency of the Slovak Republic under the grant 1/0585/19.}

\author[RUB]{Ruth Mitze}\ead{ruth.mitze@rub.de},    
\author[STU]{Michal Kvasnica}\ead{michal.kvasnica@stuba.sk},
\author[RUB]{Martin M\"onnigmann}\ead{martin.moennigmann@rub.de}               

\address[RUB]{Automatic Control and Systems Theory, Ruhr-Universit\"at Bochum, Bochum, Germany.}
\address[STU]{Institute of Information Engineering, Automation, and Mathematics, Faculty of Chemical and Food Technology, Slovak University of Technology in Bratislava, Slovakia.}


\begin{abstract}                          
This paper studies symmetric constrained linear-quadratic optimal control problems and their parametric solutions. The parametric solution of such a problem is a piecewise-affine feedback law that can be equivalently expressed as a set of active sets. We show symmetries of the optimal control problem entail symmetries of the active sets, which can be used to simplify finding the set of active sets considerably. Specifically, we improve a recently proposed method for the dynamic-programming-based enumeration of all active sets. The achieved reduction of the computational effort is illustrated with an example.
\end{abstract}

\end{frontmatter}

\section{Introduction}

The parametric solution to a constrained linear-quadratic optimal control problem (OCP), which is required for explicit linear model predictive control, for example, is a piecewise-affine feedback law \cite{Bemporad2002,Seron2003}. Every affine piece can be assigned a unique active set, so the solution can be equivalently expressed by a set of active sets. The calculation of the solution often is computationally demanding. 
Various approaches exist.
Geometric approaches exploit relations of neighboring affine pieces in the solution \cite{Bemporad2002,Tondel2003,Baotic2002}. These approaches are competitive and established in mature software toolboxes \cite{Herceg2013,Bemporad2004}.
Combinatorial approaches, also referred to as implicit enumeration techniques, calculate the set of active sets that defines the solution by taking advantage of relations of active sets \cite{Gupta2011,Feller2013,Oberdieck2017,Ahmadi2018,Herceg2015}.
Recently, the combinatorial approach was combined with dynamic programming (DP) techniques \cite{Mitze2020Aut,Monnigmann2019}.

The present paper studies constrained linear-quadratic OCPs with symmetries. Symmetries often result from symmetries of the underlying physical system and its constraints. 
Models of pendula or quadrocopters that are invariant under repeated rotations along their vertical axes by 180 or 90 degrees, respectively, are common examples (for more examples see \cite{DanielsonDiss}). 
Symmetries can concisely be described with basic group theoretic notions (see Sect.~\ref{sec:groupsAndSymmetry} for a brief summary).
These descriptions have been used to reduce the memory requirements of the explicit piecewise affine control laws that result in explicit MPC~\cite{Danielson2012,Danielson2014,Danielson2015}.
For the special case of problems with constraints that are point-symmetric to the origin, there exist pairwise symmetric active sets with similar properties \cite{Feller2013,Mitze2020IFAC}. It is an obvious question how general symmetries are reflected in the active sets.

We show that a symmetry of a constrained linear-quadratic OCP results in a symmetry in the set of active sets. Consequently, the set of all active sets can be represented by a subset thereof and the symmetric counterparts of the elements of this subset. This insight can be exploited to 
reduce the effort required to construct  all active sets.
We propose an exploration strategy for the combinatorial tree that exploits the set symmetry. The reduction achieved by these improvements is analyzed by applying this strategy to the DP approach from \cite{Mitze2020Aut}.

Section~\ref{sec:ProblemStatement} introduces constrained linear-quadratic OCPs, their symmetries, and the DP approach. Section~\ref{sec:symmetriesA} identifies symmetries in active sets, states properties of symmetric active sets, and proposes an exploration strategy for the combinatorial tree. The improved DP approach is presented in Sect.~\ref{sec:approach} 
and applied to an example in  Sect.~\ref{sec:example}. Conclusions are given in Sect.~\ref{sec:conclusion}. 

\vspace{-0.5cm}
\paragraph*{Notation.}
Consider a matrix $M\in\mathbb{R}^{a\times b}$ and an ordered set $\mathcal{M}\subseteq\{1,...,a\}$. Let $M_\mathcal{M}\in\mathbb{R}^{\vert\mathcal{M}\vert\times b}$ denote the submatrix of $M$ containing all rows indicated by $\mathcal{M}$. 
We say the sets $\mathcal{M}_1$ and $\mathcal{M}_2$ partition the set $\mathcal{M}$ if $\mathcal{M}_1\cup\mathcal{M}_2=\mathcal{M}$ and $\mathcal{M}_1\cap\mathcal{M}_2=\emptyset$. Furthermore, let $\mathcal{P}(\mathcal{M})$ refer to the power set of a set $\mathcal{M}$ and let $\min(\mathcal{M})$ and $\max(\mathcal{M})$ denote the smallest and greatest element of $\mathcal{M}$, where applicable.
Let the operators $\oplus$, $\otimes$ and $\circ$ denote the Minkowski addition, Kronecker product, and Hadamard product, respectively.
Finally, let $I^a\in\mathbb{R}^{a\times a}$, $0^{a\times b}\in\mathbb{R}^{a\times b}$, 
and $1^{a}\in\mathbb{R}^{a}$ denote the identity and zero matrices, and a column vector of ones, respectively.

\section{Problem statement and preliminaries} \label{sec:ProblemStatement}

We consider the 
constrained linear-quadratic OCP
\begin{subequations} \label{eq:OCP}
\begin{align}
\min_{U,X} \quad &\left\Vert x(N)\right\Vert_P^2+\sum_{k=0}^{N-1}\left(\left\Vert x(k)\right\Vert_Q^2+\left\Vert u(k)\right\Vert_R^2\right) \\
\textrm{s.t.} \quad 
&x(k+1)=Ax(k)+Bu(k), \: k=0,...,N-1 \label{eq:System}\\
&u(k)\in\mathcal{U}, \: k=0,...,N-1\label{eq:inputConstraints}\\
&x(k)\in\mathcal{X}, \: k=0,...,N-1\label{eq:stateConstraints}\\
&x(N)\in\mathcal{T},\label{eq:terminalConstraints}
\end{align} 
\end{subequations}
where $N\in\mathbb{N}$ is the horizon, 
inputs $u(k)\in\mathbb{R}^m$ and states $x(k)\in\mathbb{R}^n$ are collected in $U=\left(u^T(0),...,u^T(N-1)\right)^T\in\mathbb{R}^{Nm}$ and $X=\left(x^T(1),...,x^T(N)\right)^T\in\mathbb{R}^{Nn}$, respectively, 
$x(0)$ is given, 
$A\in\mathbb{R}^{n\times n}$ and $B\in\mathbb{R}^{n\times m}$ define the discrete-time time-invariant linear system, the pair $(A,B)$ is stabilizable, 
$\mathcal{U}\subset\mathbb{R}^m$ and $\mathcal{X}\subset\mathbb{R}^n$ are compact full-dimensional polytopes that contain the origin in their interiors, 
$R\in\mathbb{R}^{m\times m}$, $R\succ 0$ and $Q\in\mathbb{R}^{n\times n}$, $Q\succeq 0$ are the weighting matrices for inputs and states, respectively, the pair $(Q^{\frac{1}{2}},A)$ is detectable,
$P\in\mathbb{R}^{n\times n}$ is the optimal cost function matrix of the unconstrained infinite-horizon problem which implies $P\succ 0$, and 
$\mathcal{T}\subset\mathbb{R}^n$ is the largest possible set such that the optimal feedback for the unconstrained infinite-horizon problem stabilizes the system without violating the constraints.
Let the total number of halfspaces in \eqref{eq:OCP} and, equivalently, inequalities in \eqref{eq:qp} below be denoted $q$. Furthermore, let $\mathcal{Q}$ denote the index set $\{1,...,q\}$.

By substituting \eqref{eq:System}, \eqref{eq:OCP} can be transformed into a quadratic program (QP) of the form
\begin{align} \label{eq:qp}
\begin{split}
\min_{U} \quad &\frac{1}{2}x(0)^TYx(0)+x(0)^TFU+\frac{1}{2}U^THU\\
\textrm{s.t.} \quad &GU\leq Ex(0)+w,
\end{split}
\end{align}
with $Y\in\mathbb{R}^{n\times n}$, $F\in\mathbb{R}^{n\times Nm}$, $H\in\mathbb{R}^{Nm\times Nm}$, $G\in\mathbb{R}^{q\times Nm}$, $E\in\mathbb{R}^{q\times n}$, and $w\in\mathbb{R}^q$. We assume $w=1^q$ without restriction. The assumptions on \eqref{eq:OCP} imply $H\succ 0$ \cite{Bemporad2002}.

For any $x(0)$ such that \eqref{eq:qp} has a solution, let  $U^\star(x(0))$ denote the optimal solution to \eqref{eq:qp}.
Furthermore, let $\mathcal{A}(x(0))$ and $\mathcal{I}(x(0))$ refer to the optimal active set and inactive set, respectively, where
\begin{align*}
\mathcal{A}(x(0))&:=\{i\in\mathcal{Q}\,\vert\, G_{\{i\}}U^\star(x(0))=w_{\{i\}}+E_{\{i\}} x(0)\},\\
\mathcal{I}(x(0))&:=\mathcal{Q}\backslash\mathcal{A}(x(0)).
\end{align*}
We often drop $x(0)$. We call an active set $\mathcal{A}$ \textit{optimal} if it is optimal for some $x(0)$.

Solving \eqref{eq:qp} as a parametric program with parameter $x(0)$ results in an optimal control law $x(0)\mapsto U^\star(x(0))$ that is a continuous piecewise affine function of $x(0)$ on a polytopic partition \cite[Sect. 4.1]{Bemporad2002}. 
Let the set $\mathcal{M}_N$ collect all optimal active sets $\mathcal{A}$ such that $\mathcal{A}$ defines a full-dimensional polytope in state-space and such that $G_\mathcal{A}$ has full row rank.  Let the set $\mathcal{S}_N$ collect all optimal active sets $\mathcal{A}$ including those that define a lower-dimensional polytope (facet, vertex, etc.) and those such that $G_\mathcal{A}$ does not have full row rank ($\mathcal{S}_N\supseteq\mathcal{M}_N$).

\subsection{Group theory and symmetries of OCPs}\label{sec:groupsAndSymmetry}

We introduce some basic facts of group theory as a preparation. Details can be found in standard texts, e.g., in \cite{Lang1990}.
A set together with an operator is a group if the operator is associative, the set is closed under the operator, the set includes the identity element, and the set includes the inverse of each of its elements.
A permutation that acts on the set $\mathcal{M}$ is a bijective function that maps $\mathcal{M}$ to itself.
A permutation group is a set of permutations that satisfies the group axioms stated above. The operator is function composition.
Let a permutation group act on the set $\mathcal{M}$. The orbit $\mathcal{O}(a)$ for an $a\in\mathcal{M}$ contains the images of $a$ under all permutations that are element of the permutation group.
The orbits $\mathcal{O}(a)$ of all $a\in\mathcal{M}$ partition $\mathcal{M}$.

Following \cite[Def.~4]{Danielson2015} we call an OCP \eqref{eq:OCP} symmetric to the pair of invertible matrices $(\Theta,\Omega)$, with $\Theta\in\mathbb{R}^{n\times n}$ and $\Omega\in\mathbb{R}^{m\times m}$, if the pair $(\Theta,\Omega)$ satisfies
\begin{subequations} \label{eq:symmetricPair}
\begin{align}
\Theta A=A\Theta,&\quad \Theta B=B\Omega,\label{eq:symmetricDynamics}\\
\Theta\circ\mathcal{X}=\mathcal{X},\quad\Omega\circ\mathcal{U}&=\mathcal{U},\quad\Theta\circ\mathcal{T}=\mathcal{T},\label{eq:symmetricConstraints}\\
\Theta^TQ\Theta=Q,\quad\Omega^TR\Omega&=R,\quad\Theta^TP\Theta=P,\label{eq:symmetricCost}
\end{align}
\end{subequations}
which are the symmetries of the dynamics, constraints, and cost, respectively.
Symmetries of the dynamics and the constraints are properties of the underlying physical system. Since the weighting matrices are design parameters of the OCP~\eqref{eq:OCP}, we can often choose them to respect the desired symmetry properties. Methods for identifying all pairs $(\Theta,\Omega)$ that satisfy \eqref{eq:symmetricPair} for an OCP \eqref{eq:OCP} can be found in \cite{Danielson2014}.
Let the set $\mathcal{G}$ contain all $(\Theta,\Omega)$ that satisfy \eqref{eq:symmetricPair} for an arbitrary but fixed OCP \eqref{eq:OCP}. The set $\mathcal{G}$ is a group under pairwise matrix multiplication. The relations
\begin{subequations}
\begin{align}
&(I^n,I^m)\in\mathcal{G},\label{eq:groupIdentity}\\
&(\Theta,\Omega)\in\mathcal{G}\quad\Rightarrow\quad(\Theta^{-1},\Omega^{-1})\in\mathcal{G},\label{eq:groupInverse}\\
&(\Theta_1,\Omega_1),(\Theta_2,\Omega_2)\in\mathcal{G}\quad\Rightarrow\quad(\Theta_1\Theta_2,\Omega_1\Omega_2)\in\mathcal{G}\label{eq:groupClosure}
\end{align}
\end{subequations}
hold \cite[Thm.~2 and Proof of Prop. 1]{Danielson2015}.
A numerical example together with all pairs $(\Theta,\Omega)$ that satisfy \eqref{eq:symmetricPair} is stated in Sect.~\ref{sec:example}.

\subsection{Dynamic programming approach}\label{sec:dynamicProg}

The DP approach presented in \cite{Mitze2020Aut} solves constrained linear-quadratic OCPs by iteratively increasing the horizon from the initial horizon $N=1$ to some target horizon $N=N_{\max}$.
Without restriction, we assume the constraints in \eqref{eq:OCP} and \eqref{eq:qp} to be ordered by increasing stage
\begin{align} \label{eq:order}
\left.\begin{array}{c}
  u(0)\in\mathcal{U},\quad x(0)\in\mathcal{X},  \\ 
  \vdots  \\
  u(N-1)\in\mathcal{U},\quad x(N-1)\in\mathcal{X},  \\ 
  x(N)\in\mathcal{T}  
\end{array}\right.
\end{align}
with stages $k=0, \dots, N$. 
Let $q_\mathcal{U}$, $q_\mathcal{X}$, and $q_\mathcal{T}$ refer to the number of halfspaces that define $\mathcal{U}$, $\mathcal{X}$, and $\mathcal{T}$, respectively.
These definitions imply that $q_0= q_\mathcal{U}+q_\mathcal{X}$ and $\mathcal{Q}_0= \{1, \dots, q_0\}$ are the number and set of all constraints on $x(0)$ and $u(0)$, respectively. 
We introduce $q_0$ and $\mathcal{Q}_0$ because the number of and set of constraints on $x(0)$ and $u(0)$ play an important role in the remainder of the text.

The solution for the initial horizon $\mathcal{S}_1$ can be constructed by a minor modification of Alg.~\ref{algorithm:initHorizon} in Sect.~\ref{sec:approach}. In particular, one would omit lines 3 and 14-16, the amendment "by increasing constraints" in line 2, and the superscript "$\rm red$". 
The set $\mathcal{S}_1$ contains all elements of $\mathcal{P}(\{1,...,q_0+q_\mathcal{T}\})$ that are optimal active sets. An active set $\mathcal{A}$ is optimal if the solution for the linear program (LP)
\begin{subequations} \label{eq:FeasibilityLpWithStationarity}
\begin{align}
\min_{U,x(0),\lambda_{\mathcal{A}},s_{\mathcal{I}},t} \quad &-t \\
\textrm{s.t.} \quad & F^Tx(0) +HU+(G_{\mathcal{A}})^T\lambda_{\mathcal{A}}=0, \label{eq:LPStationarity1} \\
&t\,1^{\vert\mathcal{A}\vert}\leq \lambda_{\mathcal{A}}, \label{eq:LPStationarity2} \\ 
&G_{\mathcal{A}}U-E_{\mathcal{A}}x(0)-w_{\mathcal{A}}=0, \\
&G_{\mathcal{I}}U-E_{\mathcal{I}}x(0)-w_{\mathcal{I}}+s_{\mathcal{I}}=0, \\ 
&t\,1^{\vert\mathcal{I}\vert}\leq s_{\mathcal{I}},\, t\geq 0,
\end{align}
\end{subequations}
with 
Lagrangian multipliers $\lambda_{\mathcal{A}}$ and slack variables $s_{\mathcal{I}}$, exists \cite[Sect.\ 3.1]{Gupta2011}. If it exists, we denote the solution to \eqref{eq:FeasibilityLpWithStationarity} with $U^\star$, $x^\star(0)$, $\lambda_{\mathcal{A}}^\star$, $s_{\mathcal{I}}^\star$, $t^\star$.
An additional criterion enables to efficiently identify active sets that are not optimal: We call active sets $\mathcal{A}$ such that \eqref{eq:FeasibilityLpWithStationarity} without \eqref{eq:LPStationarity1} and \eqref{eq:LPStationarity2} has a solution (resp. has no solution) \textit{feasible} (resp. \textit{infeasible}). Feasibility is a necessary condition for optimality because the test for feasibility involves a subset of the constraints of the test for optimality.  We collect all active sets that are detected as infeasible in $\mathcal{S}_1^{\rm pruned}$. Any active set $\mathcal{A}$ such that $\mathcal{A}\supseteq\tilde{\mathcal{A}}$ for an $\tilde{\mathcal{A}}\in\mathcal{S}_1^{\rm pruned}$ is infeasible \cite[Thm.~1]{Gupta2011} and thus not optimal. Discarding active sets because they are supersets of known infeasible active sets is further referred to as pruning. Pruning has the greatest impact when the active sets are tested in the order of increasing cardinality. 
For later use, the algorithm collects active sets such that the solution to \eqref{eq:FeasibilityLpWithStationarity} is $t^\star=0$ in $\mathcal{S}_1^{\rm degen}$.

The solution $\mathcal{S}_{N+1}$ for an increased horizon is calculated with the solution $\mathcal{S}_N$. Lemma~\ref{lem:solutionIncreasedHorizon} states a basic relation of the solutions for successive horizons, which is needed below.
\begin{lem}[{\cite[Cor.~3]{Mitze2020Aut}}]\label{lem:solutionIncreasedHorizon}
	Consider an OCP \eqref{eq:OCP} with constraint order \eqref{eq:order}. 
	Assume we know $\mathcal{S}_{N}$. 	Then 
	\begin{align}\label{eq:reducedCandidatesHelper100}
	  \mathcal{S}_{{N}+1} = \mathcal{H}^{(1)}\cup\mathcal{H}^{(2)}, 
	\end{align}
	where
	\begin{subequations}
	\begin{align}
    	\mathcal{H}^{(1)}&=\left\{\mathcal{A}\in\mathcal{S}_{N}\vert\mathcal{A}\subseteq\left\{1,...,Nq_0\right\}\right\},    	\label{eq:reducedCandidatesHelper3}	\\
    	\mathcal{H}^{(2)}&\subseteq\left\{\mathcal{A}_j\cup(\mathcal{A}_l\oplus \{q_0\}) \,\vert\, 	\mathcal{A}_j\in\mathcal{P}(\mathcal{Q}_0), \, \mathcal{A}_l\in\tilde{\mathcal{S}}_N\right\},\label{eq:reducedCandidates}\\
		\tilde{\mathcal{S}}_{N}	&=\left\{\mathcal{A}\in\mathcal{S}_{N}\vert\mathcal{A}\not\subseteq\left\{1,...,(N-1)q_0\right\}\right\}. \label{eq:TildeSN}
	\end{align}
	\end{subequations}
\end{lem}
The algorithm that implements Lem.~\ref{lem:solutionIncreasedHorizon} results from Alg.~\ref{algorithm:SNp1FromSN} in Sect.~\ref{sec:approach} if the superscript "$\rm red$" is omitted.   
All elements of $S_{N+1}$ are optimal active sets by definition. Therefore, only those elements of the superset in \eqref{eq:reducedCandidates} that are optimal active sets (i.e., their solution to \eqref{eq:FeasibilityLpWithStationarity} exists) are added to $\mathcal{S}_{N+1}$. The algorithm also implements pruning to identify active sets that are not optimal more efficiently. 
The set $\mathcal{S}_{N+1}^{\rm degen}$ consists of those active sets $\mathcal{A}\in\mathcal{S}_N^{\rm degen}$ that are elements of \eqref{eq:reducedCandidatesHelper3} and those active sets in the superset in \eqref{eq:reducedCandidates} such that the solution to \eqref{eq:FeasibilityLpWithStationarity} is $t^\star=0$.

The algorithm for determining the solution $\mathcal{M}_N$ for $N=N_{\rm max}$ results from Alg.~\ref{algorithm:dynamicProgrammingApproach} in Sect.~\ref{sec:approach} if line 11 is replaced by "add $\mathcal{A}_k$ to $\mathcal{M}_N$", and if the superscript "$\rm red$" is omitted.
The algorithm iteratively increases the horizon $N$ and terminates either if $N=N_{\rm max}$ or if an $N$ such that $\mathcal{S}_N=\lim_{N\rightarrow\infty} \mathcal{S}_N$ has been found with \cite[Prop.~1]{Mitze2020Aut}. 
The set $\mathcal{S}_N$ is reduced to $\mathcal{M}_N$ by discarding those elements $\mathcal{A}$ such that $G_\mathcal{A}$ does not have full row rank and those that define lower-dimensional polytopes.
Only those elements $\mathcal{A}$ that are elements of $\mathcal{S}_N^{\rm degen}$, i.e., those that result $t^\star=0$ are tested for defining a full-dimensional polytope. The other elements are known to define a full-dimensional polytope according to \cite[Thm.~2]{Tondel2003} because  $t^\star>0$ implies strict inequality for both $G_\mathcal{I}U^\star<w_\mathcal{I}+E_\mathcal{I}x^\star(0)$ and $\lambda_\mathcal{A}^\star>0$.

\section{Symmetries in active sets}\label{sec:symmetriesA}

In this section, we identify sets of symmetric active sets and identify properties that symmetric active sets have in common. We need Thm.~\ref{thm:partnerConstraints} as a preparation. It establishes every constraint has a symmetric one for a pair $(\Theta,\Omega)\in\mathcal{G}$.
\begin{thm}\label{thm:partnerConstraints}
Consider an OCP \eqref{eq:OCP} and a pair $(\Theta,\Omega)$ such that the OCP is symmetric to the pair with \eqref{eq:symmetricPair}. Assume the OCP is reformulated as a QP \eqref{eq:qp} with $w=1^q$.
Then for every $i\in\mathcal{Q}$, there exists a $j\in\mathcal{Q}$ such that
\begin{equation}\label{eq:partnerConstraints}
\begin{split}
G_{\{i\}}&=G_{\{j\}}\left(I^N\otimes\Omega\right),\\
E_{\{i\}}&=E_{\{j\}}\Theta.
\end{split}
\end{equation}

The proof is stated in Appendix~\ref{sec:proofPartnerConstraints}.
\end{thm}

Let $\pi^{(\Theta,\Omega)}:\mathcal{Q}\rightarrow\mathcal{Q}$ refer to the function with $j=\pi^{(\Theta,\Omega)}(i)$ from \eqref{eq:partnerConstraints}.
The group properties of $\mathcal{G}$ carry over to the functions $\pi^{(\Theta,\Omega)}$ for all $(\Theta,\Omega)\in\mathcal{G}$:
\begin{thm}\label{thm:piGroup}
Consider an OCP \eqref{eq:OCP} and let $\mathcal{G}$ contain all pairs $(\Theta,\Omega)$ that satisfy \eqref{eq:symmetricPair}. 
Then, the set of all $\pi^{(\Theta,\Omega)}$, $(\Theta,\Omega)\in\mathcal{G}$ is a permutation group on $\mathcal{Q}$.

The proof is stated in Appendix~\ref{sec:proofpiGroup}.
\end{thm}

Essentially, we are interested in symmetries because they imply certain parts of our system can indistinguishably be replaced by other parts, just as exchanging the four symmetric parts of a quadrocopter by certain rotations result in indistinguishable systems. The function $\pi^{(\Theta,\Omega)}$ tells us which constraint $j=\pi^{(\Theta,\Omega)}(i)$ indistinguishably replaces constraint $i$ under such a rotation, or more generally, under a symmetry $(\Theta,\Omega)$ of our system. Just as $\pi^{(\Theta,\Omega)}$ maps \textit{single constraints} to symmetric constraints, we can map \textit{active sets} to symmetric active sets. The function $\Pi^{(\Theta,\Omega)}:\mathcal{P}(\mathcal{Q})\rightarrow\mathcal{P}(\mathcal{Q})$,
\begin{align*}
\Pi^{(\Theta,\Omega)}(\mathcal{A}):=\left\{j\in\mathcal{Q}\:\vert\: j=\pi^{(\Theta,\Omega)}(i),i\in\mathcal{A}\right\}
\end{align*}
is introduced for this purpose.
It follows from Thm.~\ref{thm:piGroup} that the set $\left\{\Pi^{(\Theta,\Omega)}\vert (\Theta,\Omega)\in\mathcal{G}\right\}$ is a permutation group on $\mathcal{P}(\mathcal{Q})$.
We call the active set $\mathcal{A}_j=\Pi^{(\Theta,\Omega)}(\mathcal{A}_i)$ the symmetric active set to $\mathcal{A}_i$ under the pair $(\Theta,\Omega)$.

The orbit $\mathcal{O}(\mathcal{A})$ of the active set $\mathcal{A}$ refers to the set of symmetric active sets to $\mathcal{A}$ under all pairs $(\Theta,\Omega)$ that satisfy \eqref{eq:symmetricPair}:
\begin{align}\label{eq:orbit}
\begin{split}
\mathcal{O}(\mathcal{A}):=\left\{\Pi^{(\Theta,\Omega)}(\mathcal{A})\:\vert\: (\Theta,\Omega)\in\mathcal{G}\:\right\}.
\end{split}
\end{align}
In other words, the orbit $\mathcal{O}(\mathcal{A})$ collects, for any active set $\mathcal{A}$, all symmetric active sets that can be generated by applying all symmetries of the given system.
It follows from the group properties of the set $\left\{\Pi^{(\Theta,\Omega)}\vert (\Theta,\Omega)\in\mathcal{G}\right\}$ that the orbits of all sets in $\mathcal{P}(\mathcal{Q})$ partition $\mathcal{P}(\mathcal{Q})$.
We drop argument $\mathcal{A}$ in the notation of an orbit and simply denote it $\mathcal{O}$. 

Theorem~\ref{thm:symmetricSameProperties} states properties that the elements of an orbit have in common.
\begin{thm}\label{thm:symmetricSameProperties}
Considering an OCP \eqref{eq:OCP}, the following statements hold:
\begin{enumerate}
\item The active sets of an orbit are either all optimal or all not optimal.
\item The polytopes defined by the active sets of an orbit are either all full-dimensional or all lower-dimensional.
\end{enumerate}

\begin{pf}
Consider two active sets $\mathcal{A}_i$ and $\mathcal{A}_j$ that are elements of the same orbit. It follows that there exists a pair $(\Theta,\Omega)\in\mathcal{G}$ such that $\mathcal{A}_j=\Pi^{(\Theta,\Omega)}(\mathcal{A}_i)$. Let the polytopes defined by $\mathcal{A}_i$ and $\mathcal{A}_j$ be denoted $\mathcal{R}_i$ and $\mathcal{R}_j$, respectively.  The pair $(\Theta,\Omega)$ causes a state-space transformation of the optimal control law with $\Theta$. Therefore, 
\begin{align}\label{eq:symmetricPolytope}
\mathcal{R}_j=\Theta\circ\mathcal{R}_i,
\end{align}
according to the proof of Cor.~1 in \cite{Danielson2015}.
Equation~\eqref{eq:symmetricPolytope} implies either both polytopes are not empty, i.e., both active sets are optimal, or both polytopes are empty, i.e., both active sets are not optimal. This proves the first claim. 
The second claim can be proved by showing that, if $\mathcal{R}_i$ is full-dimensional in state-space, then $\mathcal{R}_j$ is full-dimensional in state-space. 
Assume $\mathcal{R}_i$ is full-dimensional. Then there exist $n+1$ states $x_0,...,x_n\in\mathcal{R}_i$ such that the vectors $x_i-x_0$, $i=1,...,n$ are linearly independent, i.e.,
\begin{align}\label{eq:proofSymmetricIsFulldim1}
x_n-x_0=r_1(x_1-x_0)+...+r_{n-1}(x_{n-1}-x_0)
\end{align}
with $r_1,...,r_{n-1}\in\mathbb{R}$ has no solution.
It follows from \eqref{eq:symmetricPolytope} that $\Theta x_0,...,\Theta x_n\in\mathcal{R}_j$. The vectors $\Theta x_i-\Theta x_0$, $i=1,...,n$ are linearly independent because
\begin{align*}
\Theta x_n-\Theta x_0=&r_1(\Theta x_1-\Theta x_0)+\\
&...+r_{n-1}(\Theta x_{n-1}-\Theta x_0),
\end{align*}
equals \eqref{eq:proofSymmetricIsFulldim1} if $\Theta^{-1}$ is multiplied from the left. The inverse $\Theta^{-1}$ exists by definition of the pair $(\Theta,\Omega)$ in Sect.~\ref{sec:groupsAndSymmetry}. It follows that $\mathcal{R}_j$ is a full-dimensional state-space polytope.\hfill$\square$
\end{pf}
\end{thm}

\subsection{Exploration strategy for the combinatorial tree}

Combinatorial approaches calculate $\mathcal{M}_N$ by identifying those elements of $\mathcal{P}(\mathcal{Q})$ that are optimal active sets and define full-dimensional polytopes for a QP \eqref{eq:qp}.
The properties stated so far in Sect.~\ref{sec:symmetriesA} now allow simplifying the algorithms that exist for this purpose, because they show how to exploit the symmetries of an OCP when operating with its active sets. It suffices to analyze only one active set of a set of symmetry-related active sets. More precisely,
it is sufficient to test one active set and its polytope of each orbit for optimality and full-dimensionality to get a result for all active sets in the orbit (Thm.~\ref{thm:symmetricSameProperties}).
In the remainder of the section, we propose an exploration strategy for the active sets of $\mathcal{P}(\mathcal{Q})$ to efficiently identify those active sets that do not need to be tested.

Let a combinatorial tree (Fig.~\ref{fig:combTree}) represent all elements of $\mathcal{P}(\mathcal{Q})$. This rooted tree was first introduced in \cite[Def.~2.2]{Rymon1992} as enumeration set tree and is used in many combinatorial approaches \cite{Gupta2011,Feller2013,Ahmadi2018}, where it is referred to as active set nodal tree, active set tree, and combinatorial search tree.
Let the set $\mathcal{D}(\mathcal{A})$ contain all active sets in the subtree of $\mathcal{A}$, i.e., all descendants of $\mathcal{A}$ and $\mathcal{A}$ itself. Formally,
\begin{align}\label{eq:subtree}
\mathcal{D}(\mathcal{A}):=\left\{\mathcal{A}\cup\hat{\mathcal{A}}\,\vert\,\hat{\mathcal{A}}\in\mathcal{P}(\{\max(\mathcal{A})+1,...,q\})\right\}.
\end{align}
Furthermore, active sets with identical cardinalities appear in the same level of the tree, and each level is ordered such that those active sets that contain the lowest constraint indices appear leftmost without restriction. 

The elements of an orbit have identical cardinalities, see \eqref{eq:orbit}, and therefore appear on the same level of the combinatorial tree. 
Let the \textit{primary} active set of an orbit $\mathcal{O}$ be the element of the orbit that contains the lowest constraint indices, i.e.,
\begin{align}\label{eq:primaryA}
\mathcal{A}\in\mathcal{O}:\,\min(\mathcal{A}\backslash\mathring{\mathcal{A}})<\min(\mathring{\mathcal{A}}\backslash\mathcal{A})\,\forall\,\mathring{\mathcal{A}}\in\mathcal{O}\backslash\mathcal{A}.
\end{align}
It follows that the primary active set is the element of the orbit that is located leftmost in the level of the combinatorial tree.

Assuming the elements of the combinatorial tree are explored in the order of increasing cardinality and from left to right in each level, then the primary element of each orbit is reached first. This implies the non-primary elements need not be tested. The following Lemma can be used to efficiently identify active sets that are non-primary.
\begin{thm}\label{thm:nonprimarySubtree}
Consider an OCP \eqref{eq:OCP} and all pairs $(\Theta,\Omega)$ that satisfy \eqref{eq:symmetricPair}. If $\mathcal{A}$ is an active set that is non-primary, then all elements of $\mathcal{D}(\mathcal{A})$ are non-primary.

\begin{pf}
Consider an active set $\mathcal{A}_i$ that is non-primary. Let the primary active set in its orbit be $\mathcal{A}_p$. This implies there exists a pair $(\Theta_p,\Omega_p)\in\mathcal{G}$ such that $\Pi^{(\Theta_p,\Omega_p)}(\mathcal{A}_i)=\mathcal{A}_p$.
Furthermore, by definition of the primary active set \eqref{eq:primaryA}
\begin{align}\label{eq:proofNonprimarySubtreeHelper1}
\min(\mathcal{A}_p\backslash\mathcal{A}_i)<\min(\mathcal{A}_i\backslash\mathcal{A}_p).
\end{align}
Consider an arbitrary active set $\mathcal{A}_j\in\mathcal{D}(\mathcal{A}_i)$. It follows from \eqref{eq:subtree} that there exists an $\hat{\mathcal{A}}_i$ such that
\begin{align}\label{eq:proofNonprimarySubtreeHelper2}
\mathcal{A}_j=\mathcal{A}_i\cup\hat{\mathcal{A}}_i,\,\hat{\mathcal{A}}_i\in\mathcal{P}(\{\max(\mathcal{A}_i)+1,...,q\}).
\end{align}
Now let $\Pi^{(\Theta_p,\Omega_p)}(\mathcal{A}_j)=\mathcal{A}_k$. Then, $\mathcal{A}_k$ and $\mathcal{A}_j$ are elements of the same orbit. $\mathcal{A}_j\supseteq\mathcal{A}_i$ according to \eqref{eq:proofNonprimarySubtreeHelper2} implies $\Pi^{(\Theta_p,\Omega_p)}(\mathcal{A}_j)\supseteq\Pi^{(\Theta_p,\Omega_p)}(\mathcal{A}_i)$ and therefore,
\begin{align}\label{eq:proofNonprimarySubtreeHelper3}
\mathcal{A}_k\supseteq\mathcal{A}_p. 
\end{align}
It remains to show that $\mathcal{A}_k$ contains lower indices than $\mathcal{A}_j$, i.e.,
\begin{align}\label{eq:proofNonprimarySubtreeHelper4}
\min(\mathcal{A}_k\backslash\mathcal{A}_j)<\min(\mathcal{A}_j\backslash\mathcal{A}_k),
\end{align}
and hence, $\mathcal{A}_j$ is a non-primary active set. Inserting \eqref{eq:proofNonprimarySubtreeHelper2} and \eqref{eq:proofNonprimarySubtreeHelper3} in \eqref{eq:proofNonprimarySubtreeHelper4} gives
\begin{align}\label{eq:proofNonprimarySubtreeHelper5}
\begin{split}
\min&((\mathcal{A}_p\cup\hat{\mathcal{A}}_p)\backslash(\mathcal{A}_i\cup\hat{\mathcal{A}}_i))\\
&<\min((\mathcal{A}_i\cup\hat{\mathcal{A}}_i)\backslash(\mathcal{A}_p\cup\hat{\mathcal{A}}_p)),
\end{split}
\end{align}
with $\hat{\mathcal{A}}_i\in\mathcal{P}(\{\max(\mathcal{A}_i)+1,...,q\})$ and $\hat{\mathcal{A}}_p\in\mathcal{P}(\{1,...,q\})$.
Since $\mathcal{A}_i\neq\mathcal{A}_p$, \eqref{eq:proofNonprimarySubtreeHelper5} is equivalent to
\begin{align}\label{eq:proofNonprimarySubtreeHelper6}
\min((\mathcal{A}_p\cup\hat{\mathcal{A}}_p)\backslash\mathcal{A}_i)<\min(\mathcal{A}_i\backslash(\mathcal{A}_p\cup\hat{\mathcal{A}}_p)).
\end{align}
The expression in \eqref{eq:proofNonprimarySubtreeHelper6} is true for any $\hat{\mathcal{A}}_p$ because \eqref{eq:proofNonprimarySubtreeHelper1} holds. \hfill$\square$
\end{pf}
\end{thm}

\begin{exmp}\label{exmp:illustrationCombinatorialTree}
We illustrate Thm.~\ref{thm:nonprimarySubtree} for an artificial symmetric problem with four constraints. The combinatorial tree shown in Fig.~\ref{fig:combTree} contains all elements of $\mathcal{P}(\mathcal{Q})$, $\mathcal{Q}=\{1,...,4\}$.
\begin{figure}[tbh]
    \begin{center}
	\includegraphics[trim=0px 0 0px 0, clip, width=.45\textwidth]{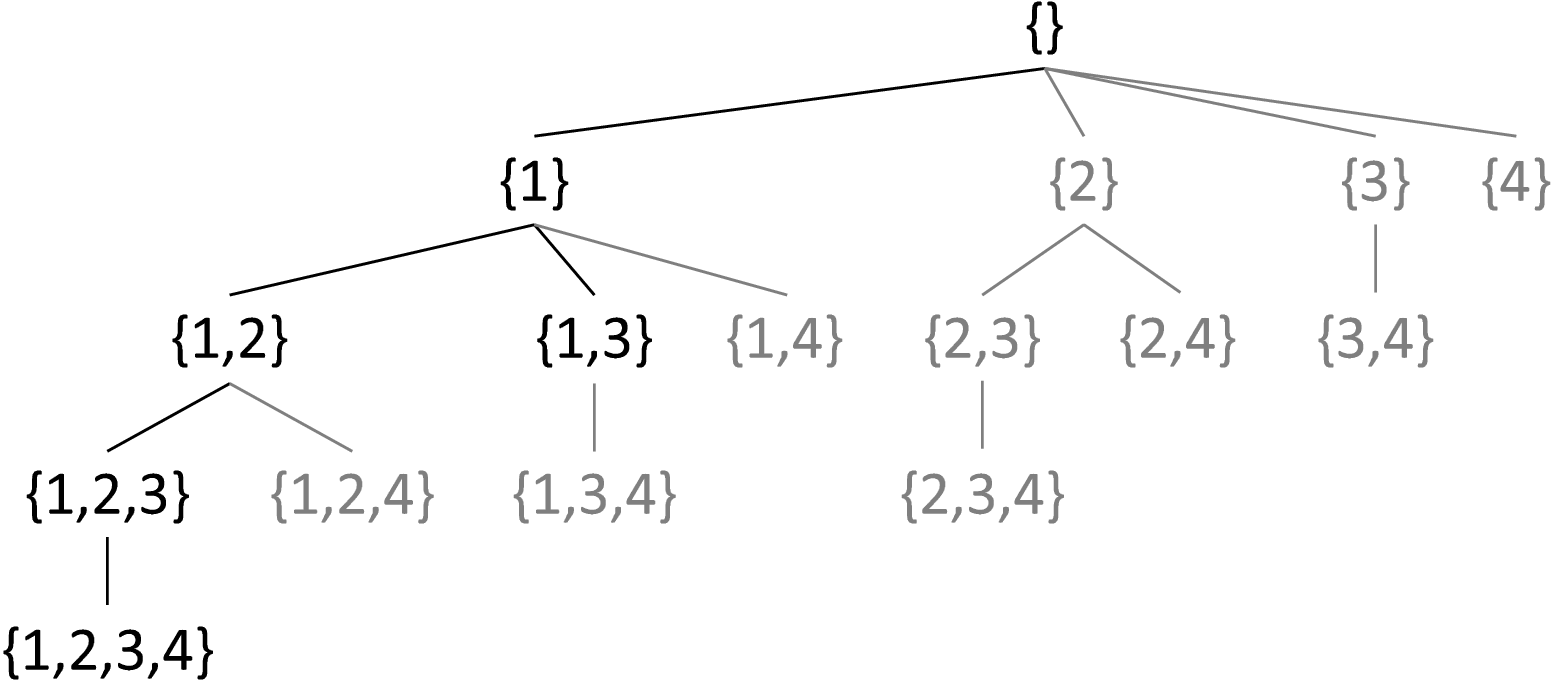}
    \caption{Combinatorial tree for Example~\ref{exmp:illustrationCombinatorialTree}. Primary active sets are shown in black, active sets that are non-primary are shown in grey.}
    \label{fig:combTree}
    \end{center}
\end{figure}
Assume there are four symmetries $(\Theta_i,\Omega_i)\in\mathcal{G}$, $i=1,...,4$ such that the functions $\pi^{(\Theta_i,\Omega_i)}$, $i=1,...,4$ are in two-line notation
\begin{align*}
\begin{split}
\pi^{(\Theta_1,\Omega_1)}&=\left(\begin{matrix}1&2&3&4\\1&2&3&4\end{matrix}\right),\quad \pi^{(\Theta_2,\Omega_2)}=\left(\begin{matrix}1&2&3&4\\2&3&4&1\end{matrix}\right),\\
\pi^{(\Theta_3,\Omega_3)}&=\left(\begin{matrix}1&2&3&4\\3&4&1&2\end{matrix}\right),\quad
\pi^{(\Theta_4,\Omega_4)}=\left(\begin{matrix}1&2&3&4\\4&1&2&3\end{matrix}\right),
\end{split}
\end{align*}
where the first row lists all elements of $\mathcal{Q}$, and for each of those elements, the second row lists the image under the permutation $\pi^{(\Theta_i,\Omega_i)}$.
The permutations are such that the set of all permutations $\pi^{(\Theta_i,\Omega_i)}$, $i=1,...,4$ is a permutation group acting on $\mathcal{Q}$.
We develop the orbit for the active set $\{1,2\}$ as an example. The orbit contains the active sets $\Pi^{(\Theta_i,\Omega_i)}(\{1,2\})$ for $i=1,...,4$, where  $\Pi^{(\Theta_i,\Omega_i)}$ maps all elements of $\{1,2\}$ with $\pi^{(\Theta_i,\Omega_i)}$. It results
\begin{align*}
\Pi^{(\Theta_1,\Omega_1)}(\{1,2\})=\{1,2\},\;\Pi^{(\Theta_2,\Omega_2)}(\{1,2\})=\{2,3\},\\
\Pi^{(\Theta_3,\Omega_3)}(\{1,2\})=\{3,4\},\;\Pi^{(\Theta_4,\Omega_4)}(\{1,2\})=\{1,4\},
\end{align*}
and therefore, 
\begin{align*}
\mathcal{O}_1=\left\{\textbf{\{1,2\}},\{2,3\},\{3,4\},\{1,4\}\right\}.
\end{align*}
The orbits for the remaining active sets in $\mathcal{P}(\mathcal{Q})$ are
\begin{align*}
\begin{split}
&\mathcal{O}_2=\left\{\textbf{\{\}}\right\},\:\mathcal{O}_3=\left\{\textbf{\{1\}},\{2\},\{3\},\{4\}\right\},\\
&\mathcal{O}_4=\left\{\textbf{\{1,3\}},\{2,4\}\right\},\\
&\mathcal{O}_5=\left\{\textbf{\{1,2,3\}},\{2,3,4\},\{1,3,4\},\{1,2,4\}\right\},\\
&\mathcal{O}_6=\left\{\textbf{\{1,2,3,4\}}\right\}.
\end{split}
\end{align*}
The set of all orbits partitions $\mathcal{P}(\mathcal{Q})$.
We marked the primary active set of each orbit in bold letters. It can be seen in Fig.~\ref{fig:combTree} that the subtrees of non-primary active sets only contain non-primary active sets as is stated in Thm.~\ref{thm:nonprimarySubtree}.
\end{exmp}

\section{Dynamic programming approach for symmetric OCP}\label{sec:approach}

We use the results from Sect.~\ref{sec:symmetriesA} to improve the DP approach from \cite{Mitze2020Aut}. 
The improvement of the algorithm is achieved by calculating only one active set of each orbit in the solution. 
The symmetric parts of the solution can then be constructed with a low computational effort.
Let $\mathcal{S}_N^{\rm red}$ denote the reduced solution set containing one active set of each orbit in
$\mathcal{A}\in\mathcal{S}_N$. 

\begin{algorithm2e}[h]
	\textbf{Initialization:} set $\mathcal{S}_1^{\rm red}=\emptyset$, $\mathcal{S}_1^{\rm degen}=\emptyset$, $S_1^{\rm pruned}= \emptyset$ and $\mathcal{S}_1^{\rm nonprim}= \emptyset$\\
	\For{every $\mathcal{A}_i\in\mathcal{P}(\{1,...,q_0+q_\mathcal{T}\})$ by incr. cardinality and by increasing constraints}{
		\If{$\mathcal{A}_i\notin\mathcal{D}(\bar{\mathcal{A}})$ for all $\bar{\mathcal{A}}\in\mathcal{S}_1^{\rm nonprim}$}
		{
			\If{$\mathcal{A}_i\not\supseteq\tilde{\mathcal{A}}$ for all $\tilde{\mathcal{A}}\in\mathcal{S}_1^{\rm pruned}$}{
				solve \eqref{eq:FeasibilityLpWithStationarity} for QP with horizon $1$\\
				\If{solution exists}{
					add $\mathcal{A}_i$ to $\mathcal{S}_1^{\rm red}$\\
					\If{solution $t^\star=0$}{
						add $\mathcal{A}_i$ to $\mathcal{S}_1^{\rm degen}$
					}
				}
				\Else{
					solve \eqref{eq:FeasibilityLpWithStationarity} without \eqref{eq:LPStationarity1} and \eqref{eq:LPStationarity2} for QP with horizon $1$\\
					\If{no solution exists}{
						add $\mathcal{A}_i$ to $S_1^{\rm pruned}$
					}
				}
			}
			\For{every $\mathcal{A}_j\in\mathcal{O}(\mathcal{A}_i)\backslash\mathcal{A}_i$}
			{
				\If{$\mathcal{A}_j\notin\mathcal{D}(\bar{\mathcal{A}})$ for all $\bar{\mathcal{A}}\in\mathcal{S}_1^{\rm nonprim}$}
				{
					add $\mathcal{A}_j$ to $\mathcal{S}_1^{\rm nonprim}$
				}
			}
		}
	}
	\textbf{Output:} $\mathcal{S}_1^{\rm red}$, $\mathcal{S}_1^{\rm degen}$
	\caption{Determination of $\mathcal{S}_1^{\rm red}$
	\label{algorithm:initHorizon}}
	\end{algorithm2e}

The algorithm for determining the reduced solution $\mathcal{S}_1^{\rm red}$ for the initial horizon is stated in Alg.~\ref{algorithm:initHorizon}. The set $\mathcal{S}_1^{\rm red}$ consists of those elements of $\mathcal{P}(\{1,...,q_0+q_\mathcal{T}\})$ that are primary and optimal.
The algorithm collects non-primary active sets in the set $\mathcal{S}_1^{\rm nonprim}$ which is initialized as empty (line 1). Any active set $\mathcal{A}$ such that $\mathcal{A}\in\mathcal{D}(\bar{\mathcal{A}})$ for an $\bar{\mathcal{A}}\in\mathcal{S}_1^{\rm nonprim}$ is non-primary with Thm.~\ref{thm:nonprimarySubtree} and thus discarded (line 3).
All active sets in $\mathcal{P}(\{1,...,q_0+q_\mathcal{T}\})$ are processed in the order of increasing cardinality and increasing constraint indices (line 2). When transferred to a combinatorial tree as in Fig.~\ref{fig:combTree}, this strategy proceeds from top to bottom, and in each level from left to right. This way, pruning has the most impact and the primary active sets of each orbit are processed first. The elements of the orbit of a primary active set, except the primary active set itself, are non-primary by definition and added to $\mathcal{S}_1^{\rm nonprim}$ if not already elements of a set $\mathcal{D}(\bar{\mathcal{A}})$, $\bar{\mathcal{A}}\in\mathcal{S}_1^{\rm nonprim}$ (lines 14-16). 
The remaining lines in the algorithm (lines 4-13) are unchanged from the DP approach.

Corollary~\ref{cor:redSolutionIncreasedHorizon} is used in Alg.~\ref{algorithm:SNp1FromSN} to construct the reduced solution $\mathcal{S}_{N+1}^{\rm red}$ for an increased horizon from the reduced solution $\mathcal{S}_N^{\rm red}$ for the current horizon. It results from modifying \cite[Cor.~5]{Mitze2020IFAC} to the general symmetries considered in this paper.
\begin{cor}\label{cor:redSolutionIncreasedHorizon}
Consider an OCP \eqref{eq:OCP} with constraint order \eqref{eq:order}. Assume we know $\mathcal{S}_{N}^{\rm red}$. Then 
	\begin{align*}
	  \mathcal{S}_{{N}+1}^{\rm red} = \mathcal{H}^{(1),\rm red}\cup\mathcal{H}^{(2),\rm red}, 
	\end{align*}
	where
	\begin{subequations}
	\begin{align}
    	\mathcal{H}^{(1),\rm red}&=\left\{\mathcal{A}\in\mathcal{S}_{N}^{\rm red}\vert\mathcal{A}\subseteq\left\{1,...,Nq_0\right\}\right\}, \\
    	\begin{split}
    	\mathcal{H}^{(2),\rm red}&\subseteq\left\{\mathcal{A}_j\cup(\mathcal{A}_l\oplus \{q_0\}) \,\vert\vphantom{\tilde{\mathcal{S}}_N^{\rm red}}\right.\\
    &\quad\quad	\left.\mathcal{A}_j\in\mathcal{P}(\mathcal{Q}_0),\mathcal{A}_l\in\tilde{\mathcal{S}}_N^{\rm red}\right\},\label{eq:redReducedCandidates} 
    	\end{split}\\
		\tilde{\mathcal{S}}_{N}^{\rm red}	&=\left\{\mathcal{A}\in\mathcal{S}_{N}^{\rm red}\vert\mathcal{A}\not\subseteq\left\{1,...,(N-1)q_0\right\}\right\}.
	\end{align}
	\end{subequations}
	
	\begin{pf}
	According to \eqref{eq:reducedCandidatesHelper100} $\mathcal{S}_{N+1}=\mathcal{H}^{(1)}\cup\mathcal{H}^{(2)}$, where $\mathcal{H}^{(1)}$ and $\mathcal{H}^{(2)}$ are defined in \eqref{eq:reducedCandidatesHelper3} and \eqref{eq:reducedCandidates}, respectively. 
    The reduced solution $\mathcal{S}_{N+1}^{\rm red}$ contains only one active set of each orbit in 
    $\mathcal{S}_{N+1}$. Let $\mathcal{F}^{(1)}$ and $\mathcal{F}^{(2)}$ contain only one active set of each orbit in 
    $\mathcal{H}^{(1)}$ and $\mathcal{H}^{(2)}$, respectively. It follows that $\mathcal{S}_{N+1}^{\rm red}=\mathcal{F}^{(1)}\cup\mathcal{F}^{(2)}$. We show that $\mathcal{F}^{(1)}=\mathcal{H}^{(1),\rm red}$ and $\mathcal{F}^{(1)}=\mathcal{H}^{(2),\rm red}$.
    $\mathcal{F}^{(1)}$ is the subset of only one active set of each orbit in
	\begin{align*}
	\mathcal{H}^{(1)}=\{\mathcal{A}\vert\mathcal{A}\subseteq\{1,...,Nq_0\},\mathcal{A}\in\mathcal{S}_N\}.
\end{align*}
Considering only one active set of each orbit in 
$\mathcal{H}^{(1)}$ can be achieved by replacing $\mathcal{S}_N$ by $\mathcal{S}_N^{\rm red}$,
\begin{align*}
\mathcal{F}^{(1)}&=\{\mathcal{A}\vert\mathcal{A}\subseteq\{1,...,Nq_0\},\mathcal{A}\in\mathcal{S}_N^{\rm red}\}\\
&=\mathcal{H}^{(1),\rm red}.
\end{align*}
In the same way, the set $\mathcal{F}^{(2)}$ results from replacing $\mathcal{S}_N$ by $\mathcal{S}_N^{\rm red}$ in
	\begin{align*}
	\mathcal{H}^{(2)}\subseteq\{\mathcal{A}\vert&\mathcal{A}=\mathcal{A}_j\cup(\mathcal{A}_l\oplus\{q_0\}),\mathcal{A}_j\in\mathcal{P}(\mathcal{Q}_0),\\
	&\mathcal{A}_l\not\subseteq\{1,...,(N-1)q_0\},\mathcal{A}_l\in\mathcal{S}_N\},
\end{align*}
which yields
	\begin{align*}
	\mathcal{F}^{(2)}\subseteq\{\mathcal{A}\vert&\mathcal{A}=\mathcal{A}_j\cup(\mathcal{A}_l\oplus\{q_0\}),\mathcal{A}_j\in\mathcal{P}(\mathcal{Q}_0),\\
	&\mathcal{A}_l\not\subseteq\{1,...,(N-1)q_0\},\mathcal{A}_l\in\mathcal{S}_N^{\rm red}\}.
\end{align*}
With $\tilde{\mathcal{S}}_{N}^{\rm red}=\left\{\mathcal{A}\in\mathcal{S}_{N}^{\rm red}\vert\mathcal{A}\not\subseteq\left\{1,...,(N-1)q_0\right\}\right\}$ this results in
	\begin{align*}
	\mathcal{F}^{(2)}&\subseteq\{\mathcal{A}\vert\mathcal{A}=\mathcal{A}_j\cup(\mathcal{A}_l\oplus\{q_0\}),\mathcal{A}_j\in\mathcal{P}(\mathcal{Q}_0),\\
	&\quad\quad\quad\mathcal{A}_l\in\tilde{\mathcal{S}}_N^{\rm red}\},
\end{align*}
and hence, $\mathcal{F}^{(2)}=\mathcal{H}^{(2),\rm red}$.
\hfill$\square$
	\end{pf}
\end{cor}
\begin{algorithm2e}[h]
	\SetKwInOut{Initialization}{Initialization}
	\SetKwInOut{Input}{Input}
	\SetKwInOut{Output}{Output}
	\textbf{Input:} $\mathcal{S}_N^{\rm red}$, $\mathcal{S}_N^{\rm degen}$\\
	\textbf{Initialization:} set $\mathcal{S}_{N+1}^{\rm red}=\emptyset$, $\mathcal{S}_{N+1}^{\rm degen}$ and $\mathcal{S}_{N+1}^{\rm pruned}= \emptyset$\\
	\For{every $\mathcal{A}_l\in\mathcal{S}_{N}^{\rm red}$}{
		\If{$\mathcal{A}_l\subseteq\{1,...,Nq_0\}$}{
			add $\mathcal{A}_l$ to $\mathcal{S}_{N+1}^{\rm red}$\\
			\If{$\mathcal{A}_l\in\mathcal{S}_N^{\rm degen}$}{
				add $\mathcal{A}_l$ to $\mathcal{S}_{N+1}^{\rm degen}$
			}
		}
		\If{$\mathcal{A}_l\not\subseteq\{1,...,(N-1)q_0\}$}{
			\For{every $\mathcal{A}_i=\mathcal{A}_j\cup(\mathcal{A}_l\oplus\{q_0\})$ with $\mathcal{A}_j\in\mathcal{P}(\mathcal{Q}_0)$  by increasing cardinality}{
				\If{$\mathcal{A}_i\not\supseteq\tilde{\mathcal{A}}$ for all $\tilde{\mathcal{A}}\in\mathcal{S}_{N+1}^{\rm pruned}$}{
					solve \eqref{eq:FeasibilityLpWithStationarity} for QP for horizon $N+1$\\
					\If{solution exists}{
						add $\mathcal{A}_i$ to $\mathcal{S}_{N+1}^{\rm red}$\\
						\If{solution $t^\star=0$}{
							add $\mathcal{A}_i$ to $\mathcal{S}_{N+1}^{\rm degen}$
						}
					}
					\Else{
						solve \eqref{eq:FeasibilityLpWithStationarity} without \eqref{eq:LPStationarity1} and \eqref{eq:LPStationarity2} for QP for horizon $N+1$\\
						\If{no solution exists}{
							add $\mathcal{A}_i$ to $\mathcal{S}_{N+1}^{\rm pruned}$ 
						}
					}
				}
			}
		}
	}
	\textbf{Output:} $\mathcal{S}_{N+1}^{\rm red}$, $\mathcal{S}_{N+1}^{\rm degen}$
	\caption{Determination of $\mathcal{S}_{N+1}^{\rm red}$ from $\mathcal{S}_N^{\rm red}$
	\label{algorithm:SNp1FromSN}}
	\end{algorithm2e}
The relation of the reduced solutions for horizon $N$ to those for horizon $N+1$ in Cor.~\ref{cor:redSolutionIncreasedHorizon} is similar to the corresponding relation in Lem.~\ref{lem:solutionIncreasedHorizon}. Therefore, the procedure for determining the reduced solution for the increased horizon in Alg.~\ref{algorithm:SNp1FromSN} is similar to the DP approach, but requires, as inputs and outputs, the reduced solutions $\mathcal{S}_N^{\rm red}$ and $\mathcal{S}_{N+1}^{\rm red}$, respectively, instead of the solutions $\mathcal{S}_N$ and $\mathcal{S}_{N+1}$. 

The overall algorithm is stated in Alg.~\ref{algorithm:dynamicProgrammingApproach}. It proceeds analogously to the DP approach, but adds the whole orbit to $\mathcal{M}_N$ whenever an active set is detected to be part of $\mathcal{M}_N$ (lines 12, 14). This is because if an active set is optimal and defines a full-dimensional polytope, the same holds for the elements of its orbit (Thm.~\ref{thm:symmetricSameProperties}).
\begin{algorithm2e}[h]
	\textbf{Input:} $\mathcal{S}_1^{\rm red}$, $\mathcal{S}_1^{\rm degen}$ (from Alg.~\ref{algorithm:initHorizon}), $N_{\max}\ge 1$\\
	\textbf{Initialization:} set $\mathcal{M}_N=\emptyset$\\
	\For{$N=1$ to $N_{\max}-1$}{
		determine $\mathcal{S}_{N+1}^{\rm red}$ and $\mathcal{S}_{N+1}^{\rm degen}$ with Alg.~\ref{algorithm:SNp1FromSN}\\
		\If{$\mathcal{S}_{N+1}^{\rm red}\subseteq\mathcal{P}(\{1,...,Nq_0\})$}{
			break 
		}
	}
	$N=N+1$\\
	\For{every $\mathcal{A}_k\in\mathcal{S}_N^{\rm red}$}{
		\If{$\text{rowrank}(G_{\mathcal{A}_k})=\vert\mathcal{A}_k\vert$}{
			\If{$\mathcal{A}_k\not\in\mathcal{S}_N^{\rm degen}$ or polytope def.\ by $\mathcal{A}_k$ full-dimensional}{
				add all elements of $\mathcal{O}(\mathcal{A}_k)$ to $\mathcal{M}_N$
			}
		}
	}
	\textbf{Output:} $\mathcal{M}_N$
	\caption{Dynamic programming approach to solving symmetric constrained linear-quadratic OCPs
	\label{algorithm:dynamicProgrammingApproach}}
	\end{algorithm2e}
	
\section{Example} \label{sec:example}

We introduce the following example.
\begin{exmp}\label{exmp:Danielson2012Exmp1}
Consider the system \cite[Example 1]{Danielson2012}
\begin{align*}
x(k+1)=\left(\begin{matrix}2&1\\-1&2\end{matrix}\right) x(k)+\left(\begin{matrix}1&0\\0&1\end{matrix}\right)u(k),
\end{align*}
with input and state constraints $-1\le u_i(k)\leq 1$, $i= 1, 2$ and $-1\le x_i(k)\leq 1$, $i= 1, 2$ respectively, 
and cost function matrices $Q=I^2$, $R=5000\cdot I^2$. The terminal cost $P$ and terminal set $\mathcal{T}$ are as described in Sect.~\ref{sec:ProblemStatement}.
\end{exmp}
Four pairs $(\Theta_i,\Omega_i)\in\mathcal{G}$, $i=1,...,4$ satisfy \eqref{eq:symmetricPair} for Example \ref{exmp:Danielson2012Exmp1},
\begin{align*}
&(\Theta_1,\Omega_1)=\left(I^2,I^2\right)\!, (\Theta_2,\Omega_2)=\left(\left(\begin{matrix}0&-1\\1&0\end{matrix}\right),\left(\begin{matrix}0&-1\\1&0\end{matrix}\right)\right)\!,\\
&(\Theta_3,\Omega_3)=\left(-I^2,-I^2\right)\!, (\Theta_4,\Omega_4)=(\Theta_2^{-1},\Omega_2^{-1}).
\end{align*}

\begin{figure}[tbh]
   \begin{center}
   \begin{subfigure}[b]{0.23\textwidth}
       \begin{center}
        \includegraphics[trim=0px 0 0px 0, clip, width=\textwidth]{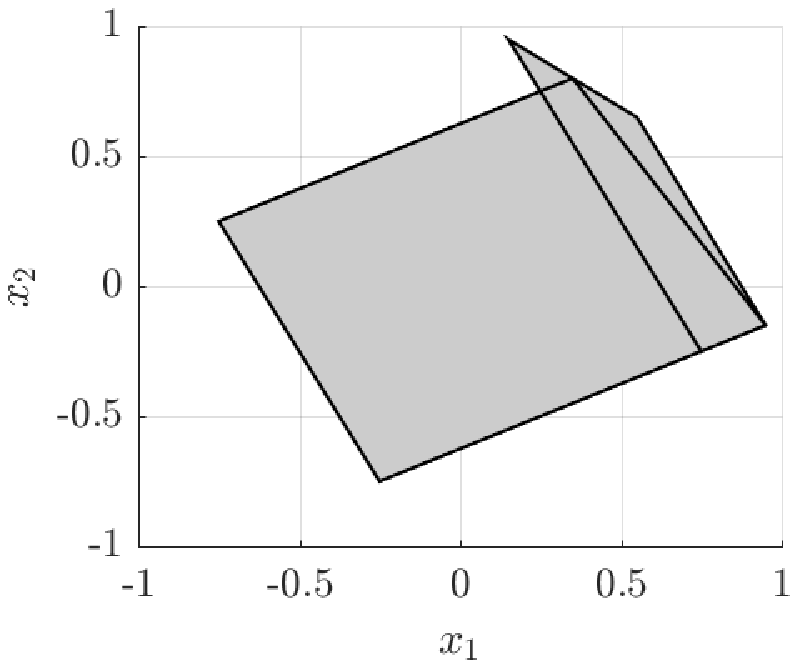}
        \caption{red. solution for ${N=1}$}
        \label{fig:Danielson2012Exmp1RedN1}
        \end{center}
	\end{subfigure}
   \begin{subfigure}[b]{0.23\textwidth}
        \begin{center}
		\includegraphics[trim=0px 0 0px 0, clip, width=\textwidth]{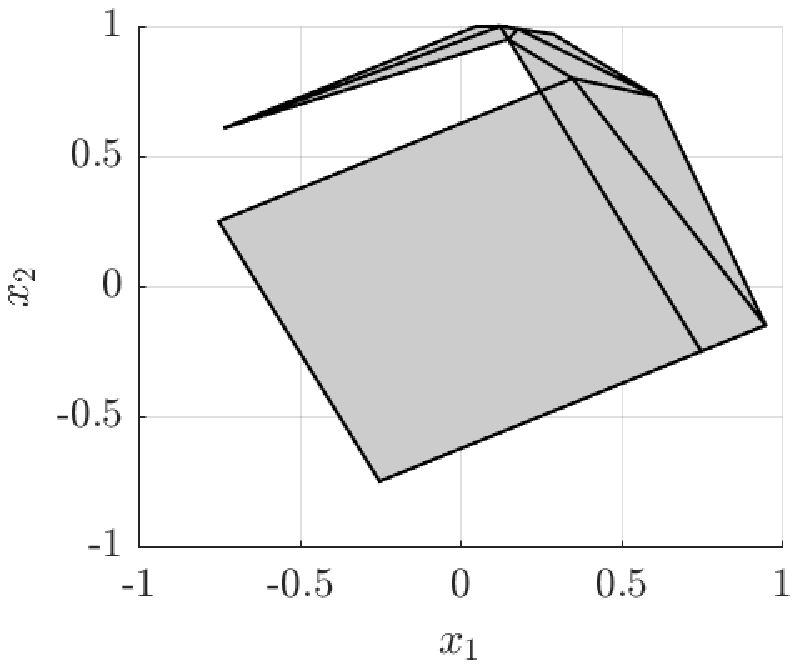}
        \caption{red. solution for ${N=2}$}
        \label{fig:Danielson2012Exmp1RedN2}
        \end{center}
	\end{subfigure}\\
	\vspace{2mm}
   \begin{subfigure}[b]{0.23\textwidth}
        \begin{center}
		\includegraphics[trim=0px 0 0px 0, clip, width=\textwidth]{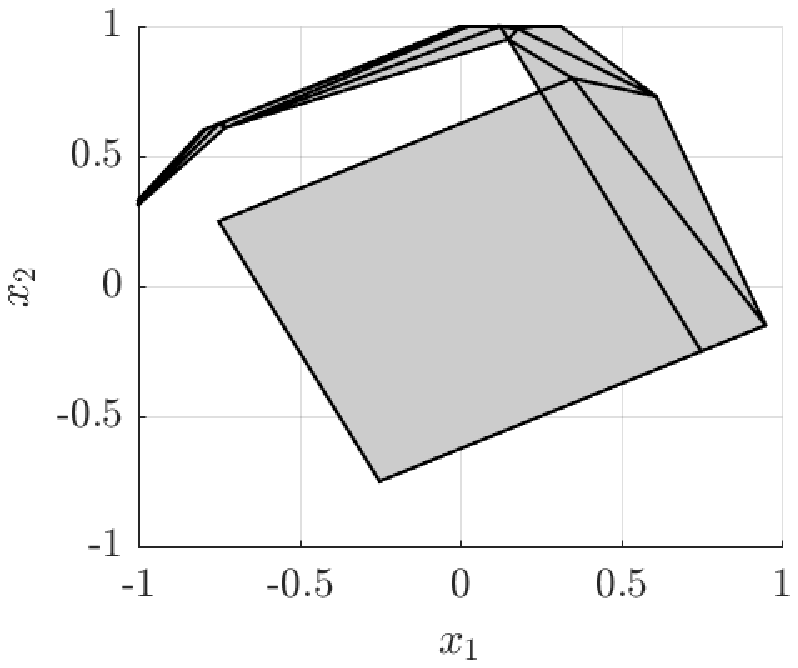}
        \caption{red. solution for ${N=5}$}
        \label{fig:Danielson2012Exmp1RedN5}
        \end{center}
	\end{subfigure}
   \begin{subfigure}[b]{0.23\textwidth}
        \begin{center}
		\includegraphics[trim=0px 0 0px 0, clip, width=\textwidth]{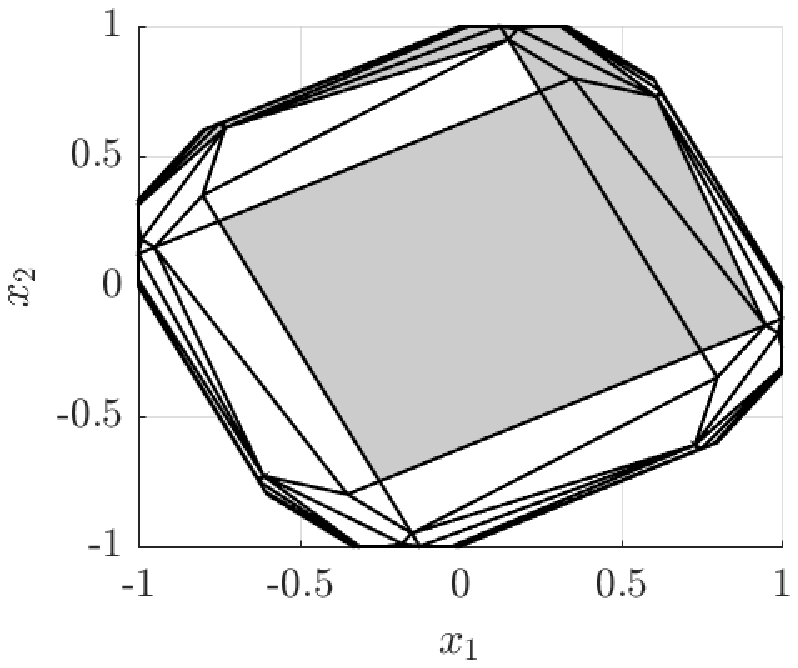}
        \caption{solution for ${N=5}$}
        \label{fig:Danielson2012Exmp1N5}
        \end{center}
	\end{subfigure}
    \caption{State-space partition for Example~\ref{exmp:Danielson2012Exmp1}. Gray polytopes are defined by active sets that are elements of the reduced solution. White polytopes are defined by the active sets in the solution that are not elements of the reduced solution in Fig.~\ref{fig:Danielson2012Exmp1N5}.}
    \label{fig:Danielson2012Exmp1Solution}
    \end{center}
\end{figure}

\vspace{-0.5cm}
\paragraph*{Construction of the solution.} We illustrate how the approach introduced in Sect.~\ref{sec:approach} constructs the solution for Example~\ref{exmp:Danielson2012Exmp1} with Fig.~\ref{fig:Danielson2012Exmp1Solution}. Figure~\ref{fig:Danielson2012Exmp1RedN1} shows the reduced solution for the initial horizon $N=1$ calculated with Alg.~\ref{algorithm:initHorizon}. The reduced solution for the increased horizon $N=2$ in Fig.~\ref{fig:Danielson2012Exmp1RedN2} is determined based on the reduced solution for $N=1$ with Alg.~\ref{algorithm:SNp1FromSN}. Repeatedly applying Alg.~\ref{algorithm:SNp1FromSN} provides the reduced solution for the target horizon, which is $N_{\max}=5$ in this example (Fig.~\ref{fig:Danielson2012Exmp1RedN5}). Figure~\ref{fig:Danielson2012Exmp1N5} shows the solution for the target horizon. It results by adding all elements of the orbits of the active sets in the reduced solution with Alg.~\ref{algorithm:dynamicProgrammingApproach}.

\vspace{-0.5cm}
\paragraph*{Computational effort.} We analyze the reduction of the computational effort that is achieved by exploiting symmetries with Table~\ref{tab:Danielson2012Exmp1CompEffort}. Table~\ref{tab:Danielson2012Exmp1CompEffort} lists the numbers of solved LPs to test for optimality with \eqref{eq:FeasibilityLpWithStationarity} and feasibility with \eqref{eq:FeasibilityLpWithStationarity} without \eqref{eq:LPStationarity1} and \eqref{eq:LPStationarity2}, and the computational times for the DP approach presented in \cite{Mitze2020Aut} and the improved DP approach presented in this paper for various horizons. The overall computational effort is dominated by the effort for solving the LPs, therefore the number of solved LPs together with the computational time suit as indicators for the computational effort. Regardless of the horizon, the number of solved LPs decreases by a factor of almost four for the approach presented in this paper. The computational time decreases by less than this factor for small horizons. In contrast, a reduction by approximately this factor is achieved for large horizons.
\begin{table}[tbh]
\caption{Computational data for Example~\ref{exmp:Danielson2012Exmp1} for various horizons $N$}
\label{tab:Danielson2012Exmp1CompEffort}
\begin{center}
\begin{tabular}{|c|c|c|c|c|}
\hline
 & \multicolumn{2}{|c|}{DP approach \cite{Mitze2020Aut}} & \multicolumn{2}{|c|}{Imp. DP approach (Sect.~\ref{sec:approach})}\\
$N$ & \# solved & comp. & \# solved & comp. \\
 & LPs & time & LPs & time \\ \hline
1 & $145$ & $1.1$\,s & $47$ ($-68$\%) & $0.8$\,s ($-27$\%)\\ \hline
3 & $2\mathord{,}917$ & $20.6$\,s & $764$ ($-74$\%) & $5.8$\,s ($-72$\%)\\ \hline
5 & $7\mathord{,}438$ & $53.2$\,s & $1\mathord{,}910$ ($-74$\%) & $14.0$\,s ($-74$\%)\\ \hline
\end{tabular}
\end{center}
\end{table}

The reduction of the number of solved LPs results for the following reasons. The set of explored active sets for the initial horizon $N=1$ is $\mathcal{P}(\{1,...,q_0+q_\mathcal{T}\})$ (line 2 in Alg.~\ref{algorithm:initHorizon}). Both approaches execute optimality tests only for a subset of these active sets. The DP approach discards those that are supersets of known infeasible active sets, the approach presented in this paper additionally discards those that are non-primary (line 3 in Alg.~\ref{algorithm:initHorizon}). It follows that the algorithm presented in this paper executes fewer optimality tests. Feasibility tests are executed for all active sets that were tested for optimality and were identified as not optimal. Thus, fewer optimality tests entail fewer feasibility tests. 
For any increased horizon $N+1$, the number of explored active sets for the DP approach depends on the cardinality of the solution $\mathcal{S}_N$. 
For the approach presented in this paper, it depends on the cardinality of the reduced solution $\mathcal{S}_N^{\rm red}$ (line 3 in Alg.~\ref{algorithm:SNp1FromSN}) which includes only one active set of each orbit in $\mathcal{S}_N$.
Thus, there result fewer explored active sets for the approach presented in this paper, and consequently, fewer optimality and feasibility tests are executed.
Fewer solved LPs result in less computational time. However, the DP approach presented in this paper requires additional computational effort, most notably for identifying and discarding non-primary active sets for the initial horizon (lines 3,14-16 in Alg.~\ref{algorithm:initHorizon}). This effort diminishes the reduction of computational time. Its relative impact is only relevant for short horizons, however.

It follows that the higher the average number of elements of the orbits, the higher the reduction of the computational effort due to the decreased number of LPs. The number of elements of an orbit increases with the number of the symmetries of the problem.
Considering a problem with $g$ symmetries, each orbit contains up to $g$ elements, see \eqref{eq:orbit}, and the number of solved LPs by the approach presented in this paper decreases to about $1/g$  of the solved LPs by the DP approach.
The computational effort required for additional computations by the proposed approach counteracts this reduction in computational effort and has the most impact for short horizons.

\section{Conclusion} \label{sec:conclusion}

In this work, it was shown that symmetries in constrained linear-quadratic OCPs entail symmetries in the set of active sets. 
Based on this, a strategy for improving combinatorial approaches for constrained linear-quadratic OCPs with symmetries was proposed. The computational effort for the improved approach is almost reduced to $1/g$ of the computational effort for the approach that does not consider symmetries, where $g$ denotes the number of symmetries of the OCP.

The results of this paper can be also used to implement the orbit controller introduced in \cite[Sect. 5]{Danielson2015} with less computational effort. 
A representation of the solution with only one active set of each orbit naturally results from the approach proposed here. Postprocessing steps for selecting representative active sets or pieces of the piecewise-affine solution \cite[Sect. 5B]{Danielson2015} are therefore no longer required.
This topic will be investigated in future.

\appendix

\section{Proof of Thm.~\ref{thm:partnerConstraints}}
\label{sec:proofPartnerConstraints}

\begin{pf}
In this proof, let $1$ denote a column vector of ones of appropriate size. $\mathcal{U}$, $\mathcal{X}$, and $\mathcal{T}$ can be defined by intersections of a finite number of halfspaces. Since $\mathcal{X}$, $\mathcal{U}$, and $\mathcal{T}$ are symmetric to the pair $(\Theta,\Omega)$ by assumption, the conditions \eqref{eq:symmetricConstraints} hold. As a consequence, there exist $T^\mathcal{U}$, $T^\mathcal{X}$, and $T^\mathcal{T}$ 
such that 
\begin{equation}\label{eq:proofSymQpHelper1}
  \begin{split}
	\mathcal{U}&=\left\{u\in\mathbb{R}^m\:\vert\:
	\left(\begin{matrix}T^\mathcal{U}\\ T^\mathcal{U}\Omega\end{matrix}\right)u\leq
	1\right\},\\
	\mathcal{X}&=\left\{x\in\mathbb{R}^n\:\vert\:
	\left(\begin{matrix}T^\mathcal{X}\\ T^\mathcal{X}\Theta\end{matrix}\right)x\leq
	1\right\},\\
	\mathcal{T}&=\left\{x\in\mathbb{R}^n\:\vert\:
	\left(\begin{matrix}T^\mathcal{T}\\ T^\mathcal{T}\Theta\end{matrix}\right)x\leq
	1\right\}.
	\end{split}
\end{equation}
We now investigate how the expressions in \eqref{eq:proofSymQpHelper1} relate to the input, state, and terminal constraints in the OCP \eqref{eq:OCP}. The input constraints \eqref{eq:inputConstraints} expressed with \eqref{eq:proofSymQpHelper1} yield
\begin{align}\label{eq:symmetriesOfU}
	&\left(\begin{matrix} \left(\begin{matrix} 0^{1\times k} & 1 & 0^{1\times (N-k-1)} \end{matrix}\right)\otimes T^\mathcal{U}\\
	\left(\begin{matrix} 0^{1\times k} & 1 & 0^{1\times (N-k-1)} \end{matrix}\right)\otimes T^\mathcal{U}\Omega \end{matrix}\right)U\leq 1
\end{align}
or equivalently
\begin{align*}
	&\left(\begin{matrix} 
	\left(\begin{matrix} 0^{1\times k} & 1 & 0^{1\times (N-k-1)} \end{matrix}\right)\otimes T^\mathcal{U} \\
	\left[\left(\begin{matrix} 0^{1\times k} & 1 & 0^{1\times (N-k-1)} \end{matrix}\right)\otimes T^\mathcal{U}\right] 
	\cdot\left(I^N\otimes\Omega\right)
	\end{matrix}\right)U\leq 1,
\end{align*}
where $k= 0, \dots, N-1$. 
As preparation for dealing with the state and terminal constraints, we express a state $x(k)$, $k>0$ with $x(0)$ and the input variables $U$ with \eqref{eq:System},
\begin{equation*}
x(k)=
\left(\begin{matrix} A^{k-1}B & A^{k-2}B & \cdots & B & 0^{n\times(N-k)n} \end{matrix}\right)U
	+A^{k} x(0).
\end{equation*}
The state constraints \eqref{eq:stateConstraints} expressed with \eqref{eq:proofSymQpHelper1} then yield
\begin{subequations}
\begin{align}
	&\left(\begin{matrix} T^\mathcal{X} \\ T^\mathcal{X}\Theta \end{matrix}\right)x(0)
	\leq 1,\label{eq:symmetriesOfX0}
	\\
	&\left(\begin{matrix} T^\mathcal{X}\left(\begin{matrix} A^{k-1}B & A^{k-2}B & \cdots & B & 0^{n\times(N-k)n} \end{matrix}\right) \\ T^\mathcal{X}\Theta\left(\begin{matrix} A^{k-1}B & A^{k-2}B & \cdots & B & 0^{n\times(N-k)n} \end{matrix}\right) \end{matrix}\right)U
	\nonumber
	\\
	&\quad\quad\quad+\left(\begin{matrix} T^\mathcal{X}A^k \\ T^\mathcal{X}\Theta A^k \end{matrix}\right)x(0)
\leq 1,\label{eq:symmetriesOfXk}
\end{align}
\end{subequations}
for $k= 1, \dots, N-1$. 
Submitting \eqref{eq:symmetricDynamics} into \eqref{eq:symmetriesOfXk} results in
\begin{equation}\label{eq:symmetriesOfX}
\begin{split}
	&\left(\begin{matrix} T^\mathcal{X}\left(\begin{matrix} A^{k-1}B & A^{k-2}B & \cdots & B & 0^{n\times(N-k)n} \end{matrix}\right) \\ T^\mathcal{X}\left(\begin{matrix} A^{k-1}B & A^{k-2}B & \cdots & B & 0^{n\times(N-k)n} \end{matrix}\right)\cdot\left(I^N\otimes\Omega\right) \end{matrix}\right)\\
	&\quad\quad\cdot U+\left(\begin{matrix} T^\mathcal{X}A^k \\ T^\mathcal{X}A^k\Theta \end{matrix}\right)x(0)
\leq 1,\: k=1,...,N-1.
\end{split}
\end{equation}
The constraints on the terminal stage \eqref{eq:terminalConstraints} expressed with \eqref{eq:proofSymQpHelper1} yield
\begin{equation*}
\begin{split}
	&\left(\begin{array}{c}T^\mathcal{T}\left(\begin{array}{c c c}A^{N-1}B & \cdots & B\end{array}\right)\\T^\mathcal{T}\Theta\left(\begin{array}{c c c}A^{N-1}B & \cdots & B\end{array}\right)\end{array}\right)U\\
	&+\left(\begin{array}{c}T^\mathcal{T}A^N\\T^\mathcal{T}\Theta A^N\end{array}\right)x(0)\leq
	1.
\end{split}
\end{equation*}
Submitting \eqref{eq:symmetricDynamics} results in 
\begin{equation}\label{eq:symmetriesOfT}
\begin{split}
	&\left(\begin{array}{c}T^\mathcal{T}\left(\begin{array}{c c c}A^{N-1}B & \cdots & B\end{array}\right)\\T^\mathcal{T}\left(\begin{array}{c c c}A^{N-1}B & \cdots & B\end{array}\right)\cdot\left(I^N\otimes\Omega\right)\end{array}\right)U\\
	&\quad\quad\quad+\left(\begin{array}{c}T^\mathcal{T}A^N\\T^\mathcal{T}A^N\Theta\end{array}\right)x(0)\leq
	1.
\end{split}
\end{equation}
The form $GU\leq Ex(0)+w$, $w=1$ can be obtained by moving all terms that depend on $x(0)$ to the right-hand side in \eqref{eq:symmetriesOfU}, \eqref{eq:symmetriesOfX0}, \eqref{eq:symmetriesOfX}, and \eqref{eq:symmetriesOfT}. $GU\leq Ex(0)+w$ inherits the alternating rows of \eqref{eq:symmetriesOfU}, \eqref{eq:symmetriesOfX0}, \eqref{eq:symmetriesOfX}, and \eqref{eq:symmetriesOfT}, so for every constraint $i$ there is a constraint $j$ such that $G_{\{i\}}=G_{\{j\}}\left(I^N\otimes\Omega\right)$ and $E_{\{i\}}=E_{\{j\}}\Theta$.
\hfill$\square$
\end{pf}

\section{Proof of Thm.~\ref{thm:piGroup}}\label{sec:proofpiGroup}

\begin{pf}
We show that the set of all $\pi^{(\Theta,\Omega)}$, $(\Theta,\Omega)\in\mathcal{G}$ satisfies the group axioms.

\textit{Associativity}: The composition of functions is associative, see, e.g., \cite[p.~28]{Lang1990}.

\textit{Identity}: Consider the pair $(I^n,I^m)$ and a constraint $i\in\mathcal{Q}$. For $j=\pi^{(I^n,I^m)}(i)$ the conditions in \eqref{eq:partnerConstraints} are
\begin{align*}
G_{\{i\}}&=G_{\{j\}}\left(I^N\otimes I^m\right)&\,&\Leftrightarrow\, &G_{\{i\}}&=G_{\{j\}},\\
E_{\{i\}}&=E_{\{j\}} I^n&\,&\Leftrightarrow\, &E_{\{i\}}&=E_{\{j\}}.
\end{align*}
This implies $\pi^{(I^n,I^m)}(i)=i$. With \eqref{eq:groupIdentity} this implies $(I^n,I^m)\in\mathcal{G}$.

\textit{Invertibility}: Consider a pair $(\Theta,\Omega)\in\mathcal{G}$ and a constraint $i\in\mathcal{Q}$. For $j=\pi^{(\Theta,\Omega)}(i)$ the conditions in \eqref{eq:partnerConstraints} are
\begin{align*}
G_{\{j\}}&=G_{\{i\}}\left(I^N\otimes\Omega^{-1}\right),\\
E_{\{j\}}&=E_{\{i\}}\Theta^{-1}.
\end{align*}
This implies $i=\pi^{(\Theta^{-1},\Omega^{-1})}(j)$. With \eqref{eq:groupInverse} this implies $(\Theta^{-1},\Omega^{-1})\in\mathcal{G}$.

\textit{Closure}: Consider two pairs $(\Theta_1,\Omega_1),(\Theta_2,\Omega_2)\in\mathcal{G}$ and a constraint $i\in\mathcal{Q}$. For the function composition $j=\pi^{(\Theta_2,\Omega_2)}\left(\pi^{(\Theta_1,\Omega_1)}(i)\right)$ the conditions in \eqref{eq:partnerConstraints} are
\begin{align*}
G_{\{i\}}&=\left(G_{\{j\}}\left(I^N\otimes\Omega_1\right)\right)\left(I^N\otimes\Omega_2\right)\\
&\Leftrightarrow\;G_{\{i\}}=G_{\{j\}}\left(I^N\otimes\left(\Omega_1\Omega_2\right)\right),\\
E_{\{i\}}&=\left(E_{\{j\}}\Theta_1\right)\Theta_2\;\Leftrightarrow\;E_{\{i\}}=E_{\{j\}}\left(\Theta_1\Theta_2\right).
\end{align*}
This implies $j=\pi^{(\Theta_1\Theta_2,\Omega_1\Omega_2)}(i)$. With \eqref{eq:groupClosure} this implies $(\Theta_1\Theta_2,\Omega_1\Omega_2)\in\mathcal{G}$.
\hfill$\square$
\end{pf}

\bibliographystyle{plain}        
\bibliography{root}           

\end{document}